\documentstyle{amsppt}
\magnification1200
\NoBlackBoxes

\pageheight{9 true in} \pagewidth{6.5 true in}

\topmatter
\title
Smooth values of the iterates of the Euler's phi-function
\endtitle
\author
Youness Lamzouri
\endauthor
\address{D{\'e}partment  de Math{\'e}matiques et Statistique,
Universit{\'e} de Montr{\'e}al, CP 6128 succ Centre-Ville,
Montr{\'e}al, QC  H3C 3J7, Canada}
\endaddress
\email{Lamzouri{\@}dms.umontreal.ca}
\endemail
\thanks
\endthanks
\abstract  Let $\phi(n)$ be the Euler-phi function, define
$\phi_0(n) = n$ and $\phi_{k+1}(n)=\phi(\phi_{k}(n))$ for all
$k\geq 0$. We will determine an asymptotic formula for the set of
integers $n$ less than $x$ for which $\phi_k(n)$ is $y$-smooth,
conditionally on a weak form of the Elliott-Halberstam conjecture.
\endabstract

\endtopmatter

\document

\head 1.Introduction \endhead

\noindent Integers without large prime factors, usually called
$smooth$ $numbers$, play a central role in several topics of
number theory. From multiplicative questions to analytic methods,
they have various and wide applications, and understanding their
behavior will have important consequences for number theoretic
algorithms, which are an important tool in cryptography.

\noindent Let $\phi(n)$ be the Euler-phi function, define
$\phi_0(n) = n$ and $\phi_{k+1}(n)=\phi(\phi_{k}(n))$ for all
$k\geq 0$. There are several interesting results on the behavior
of the functions $\phi_k$ (Erd\"os, Granville, Pomerance and Spiro
[5]). It is known that the understanding of the multiplicative
structure of the phi-function and its iterates is in some sense
equivalent to studying the behavior of the integers of the form
$p-1$ where $p$ is prime. It is also believed that the
distribution of the prime factors of such an integer behaves like
that of a random integer, in the following sense: Define
$$\Psi(x,y) = \big|\{n\leq x : p|n \implies  p\leq y \}\big|\quad
\hbox{ and }  \quad \pi(x,y) = \big|\{p\leq x : q|p-1 \implies  q\leq y \}\big| .$$
\proclaim {Conjecture 1} Fix $U\geq 1$. If $x^{1/U}\leq y\leq x$ then
$$  \frac{\pi(x,y)}{\pi(x)} \sim\frac{\Psi(x,y)}{x}\quad as \quad x\rightarrow \infty . $$
\endproclaim
Assuming this conjecture one can deduce the behavior of the function $\pi(x,y)$ from the known asymptotic formula
$$ \Psi(x,y) \sim x \rho(u) \hbox{ as } x \rightarrow \infty \hbox{ with } x= y^u $$
where $\rho(u)$ is the Dickman function, defined as the unique continuous solution of the differential-difference
equation $u\rho'(u) = -\rho(u-1)$ for $u\geq 1$, satisfying the initial condition $\rho(u)=1$ for $0\leq u\leq 1$.

Now let $P$ be a set of prime numbers and define
$$\Psi(x,P) = \big|\{n\leq x : p|n \implies  p\in P \}\big| \quad \hbox{ and } \quad \pi(x,P) =
 \big|\{p\leq x : q|p-1 \implies  q\in P \}\big| .$$ One might guess
$$   \frac{\pi(x,P)}{\pi(x)} \sim\frac{\Psi(x,P)}{x}\quad as \quad x\rightarrow \infty, \tag{1}$$
under certain conditions on the set $P$.

Granville [7] has an unpublished argument that Conjecture 1 holds
for $u=\log(x)/\log(y)$ bounded, assuming the Elliott-Halberstam
conjecture (E-H) which states that:
$$ \sum_{q\leq x^{1-\epsilon}} \max_{y\leq x} \max_{(a,q)=1} \bigg |\pi(y;q,a)-\frac{\pi(y)}{\phi(q)}
\bigg| \ll_{\epsilon,A}
\frac{x}{\log(x)^{A}}. $$
A weak version of this conjecture is the following:

\proclaim {Conjecture 2}Fix $\epsilon>0$. Then
$$ \sum\Sb d\leq x^{1-\epsilon} \endSb
\bigg |\pi(x;d,1)-\frac{\pi(x)}{\phi(d)}\bigg | =
o\left(\pi(x)\right) \quad \hbox{ as }\quad x\to \infty.$$
\endproclaim

We will prove a version of (1) assuming this Conjecture;
specifically we show the following: \proclaim {Theorem 1} Assume
Conjecture 2. If $P$ is a set of primes less than $x$ for which
 $$\sum\Sb p\notin P \\ p\leq x \endSb \frac{1}{p} \ll  1 \qquad  \hbox{then}\quad \frac{\pi(x,P)}{\pi(x)} \sim
 \prod_{p\notin P}\left(1-\frac{1}{(p-1)^2}\right)\frac{\Psi(x,P)}{x} \quad \hbox{ as } x\to \infty.$$

\endproclaim

Note that there is an extra factor in Theorem 1 compared with (1).
To see why we should expect this, let $q$ be some prime; then the
probability that a random integer $n$ is divisible by $q$ is
$1/q$. Now the probability that a random integer of the form $p-1$
(where $p$ prime) is divisible by $q$ is  $1/(q-1)$ (since $p$ is
excluded from the class $0\mod q$). The differences between the
two probabilities are negligible as $q$ increases, however this is
not true for small primes $q$, and thus we need a correction
factor (it can be removed in some special cases, see Lemma 2.1).

\noindent Define
$$ \Phi_k(x,y) = \big|\{n\leq x : p|\phi_k(n) \implies p\leq y\}\big|.$$
Using Theorem 1 we get an asymptotic of this function
conditionally on conjecture 2. \proclaim{Theorem 2}Assume
Conjecture 2. Fix $U>1$. If $y=x^{1/u}$ where $1\leq u\leq U$,
then
$$ \Phi_k(x,y) \sim x\sigma_k(u) \quad \hbox{as} \quad x\rightarrow \infty $$
where  $\sigma_k(u) = 1$ for  $u\leq 1$, and
$\displaystyle{u\sigma_{k+1}(u)= \int_0^u \sigma_{k+1}(u-t)
\sigma_{k}(t)dt}$ for $ u\geq 1$, with $\sigma_0(u) =
\rho(u)=((e+o(1))/u\log(u))^u$. Moreover, for all $k\geq 1$
$$\quad\sigma_k(u)=\left(\frac{1+o(1)}{\log_k(u)\log_{k+1}(u)}\right)^u \hbox{ and }
\log_k(u)=\log(\log(\log(...\log(u)...)))\quad k \hbox{ times.}
$$
\endproclaim

The first step in the proof uses simple combinatorics to
approximate the functions $\Phi_k(x,y)$ by $\Psi(x,P_k)$, where
$P_k$ are the sets of primes defined iteratively by
$P_{k+1}=\{p\leq x :q|p-1\implies q\in P_k\}$, with $P_0 = \{p\leq
y\}$. \proclaim {Proposition 1}
$$\Phi_k(x,y)=\Psi(x,P_k) + O\left(\frac{x(\log x)^{2k}}{y}\right).$$
\endproclaim
From the fact that $|P_k|=\pi(x,P_{k-1})$, the next step in
proving our Theorem 2 is to establish a relation between $|P|$ and
$\Psi(x,P)$ for any given set of primes $P$. This was done by
Granville and Soundararajan [8] while studying mean values of
multiplicative functions. They proved the following proposition:
\proclaim{ Proposition 2 (Proposition 1 of [8])}Let $f$ be a
multiplicative function with $|f(n)|\leq 1$ for all $n$, and
$f(n)=1$ for $n\leq y$. Let $\theta(x) = \displaystyle{\sum_{p\leq
x}\log(p)}$ and define
$$\chi(u) := \frac{1}{\theta(y^u)}\sum_{p\leq y^u}
f(p)\log(p).$$Then $\chi(t)$ is a measurable function with
$\chi(t) = 1$ for all $t\leq 1$. Let $\sigma$ be the corresponding
unique solution to the equation:
$$ u\sigma(u) = \int_0^u \sigma(u-t)\chi(t)dt\quad \hbox{for} \quad u>1 \tag{2}$$ subject to the initial
condition $\sigma(u)=1$ for $0 \leq u\leq 1$. Then
$$ \frac{1}{y^u}\sum_{n\leq y^u} f(n) = \sigma(u) + O\left(\frac{u}{\log(y)}\right).$$
\endproclaim From this result  and by partial summation we can deduce

\proclaim{Corollary 1} Fix $U>1$. Let $P$ be a set of primes less
than $x$ such that $P_0 \subseteq P$, and $f$ be a completely
multiplicative function such that $f(p)=1$ if $p\in P$ and $0$
otherwise (so that $f(n)=1$ for all $n\leq y$). For $1\leq u\leq
U$, define
$$ \chi(u):=\frac{1}{\pi(y^u)}\sum\Sb p\in P \\ p\leq y^u \endSb 1, $$
then
$$\Psi(y^u,P)=\sum_{n\leq x}f(n) \sim y^u\sigma(u) $$where $\sigma$ is the corresponding solution to (2).
\endproclaim
It remains to study (2), a delay integral equation, and to try to
estimate the solution $\sigma$ where $\chi$ is a certain
measurable function. In several interesting cases $\chi(u)$ decays
like $\left(\{1+o(1)\}/h(u)\right)^u$ where $h$ is positive and
non-decreasing. We prove the following:

\proclaim {Theorem 3} Let $\chi$ be a real measurable function for
which $\chi(t)=1$ for $0\leq t \leq 1$, and $0\leq\chi(t)\leq 1$
for $t>1$. Moreover suppose that

\text{\rm i)}  $\int_T^{\infty}\chi(t)dt=0$ for some constant $T$.
We define $T=\min \{t :\int_T^{\infty}\chi(t)dt=0\}$ to avoid
redundancy, and suppose that $T>1$.

\noindent or \text{\rm ii)}
$\chi(t)=\left(\{1+o(1)\}/h(t)\right)^t$ where $h(t)$ is
non-decreasing  and $h(t)\rightarrow \infty$ as $t\rightarrow
\infty$.

Let $\sigma$ be the corresponding solution to (2). Then
$$\sigma(u)=\exp\left((-\xi(u)+o(1)) u+\int_1^{\infty}\frac{\chi(v) e^{\xi(u) v}}{v}
dv\right),$$ where $\xi(u)$ is the unique solution to
$\displaystyle{u=\int_1^{\infty}\chi(v) e^{\xi(u) v} dv}.$
\endproclaim
Moreover we can get explicit asymptotic in a number of interesting
cases, we prove

\proclaim{Proposition 3} Let $\chi$ be a real measurable function
for which $\chi(t) = 1$ for $0 \leq t\leq 1$, and $0\leq
\chi(t)\leq 1$ for $t>1$. Suppose that
$\int_T^{\infty}\chi(t)dt=0$ for some constant $T$. We define
$T=\min \{t :\int_T^{\infty}\chi(t)dt=0\}$ to avoid redundancy,
and suppose that $T>1$. \noindent Then $\xi(u) =
\displaystyle{\frac{\log(u)}{T}}(1+o(1))$, and
$$\sigma(u)=\exp\left(-\frac{u\log(u)}{T}(1+o(1))\right).$$
\endproclaim

\proclaim {Proposition 4} Let $\chi$ be a real measurable function
for which $\chi(t) = 1$ for $0 \leq t\leq 1$, and $0\leq
\chi(t)\leq 1$ for $t>1$; and suppose that  $\chi(u)
=\left(\{1+o(1)\}/h(u)\right)^u$ where $h$ satisfies the following
conditions:

\noindent \text{\rm i)}\ $h$ is positive and non-decreasing with
$h(u)\rightarrow \infty$ as $u\rightarrow \infty$.

\noindent \text{\rm ii)}\ $h$ is continuously differentiable and
$uh'(u)/h(u)\rightarrow n$  as $u\rightarrow\infty$ for some
$0\leq n < \infty$.

We distinguish two cases: a) $0<n<\infty$ and b) $n=0$.

\noindent Then
 $$\sigma(u)= \left(\frac{1+o(1)}{h(\zeta\log(u))}\right)^u, $$where $\zeta=e/n$ in case a) and $\zeta=1$ in case b).

\endproclaim

The distinction between cases a) and b) in Proposition 4 justifies
the appearance of the constant $e$ only in the asymptotic of
$\sigma_o$ in Theorem 2.
\bigskip

\noindent {\bf Acknowledgments.}

 I sincerely thank my advisor,
Professor Andrew Granville, for all his advices and
encouragements. I would also thank the referee for his rigorous
comments and remarks.

\head 2. Proof of Theorem 2 \endhead

\proclaim {Lemma 2.1} Assume Conjecture 2. Fix $U\geq 1$. Suppose
that $P$ is a set of primes less than $x$ for which $\{p\leq y
\}\subseteq P$, where $y=x^{1/u}$ and $1 \leq u\leq U$. Then $$
\frac{\pi(x,P)}{\pi(x)}\sim \frac{\Psi(x,P)}{x}\quad as \quad
x\rightarrow \infty. $$
\endproclaim
\demo{Proof} We have that
$$ \sum \Sb p\notin P\\ p\leq x \endSb \frac{1}{p} \leq \sum \Sb  y < p\leq x \endSb \frac{1}{p}\ll
\log\left(\frac{\log(x)}{\log(y)}\right)=\log(u) \ll 1$$ and, since $1-t \geq e^{-2t}$ for $0\leq t\leq 1/2$,
then $$1\geq\prod_{p\notin P}\left(1-\frac{1}{(p-1)^2}\right)\geq\prod_{p>y}\left(1-\frac{1}{(p-1)^2}\right)\geq
\exp\left(-2\sum_{p>y} \frac{1}{(p-1)^2}\right) = 1+o(1). $$ The result follows by Theorem 1.
\enddemo

 \demo{Proof of Theorem 2} First note that the sets $P_k$
for $k\geq 0$ satisfy the conditions of Lemma 2.1. Now $
\Psi(x,P_0) = \Psi(x,y) \sim \rho(u)x \hbox{ as } x\rightarrow
\infty $. We use induction on $k$: suppose that $\Psi(x,P_k)\sim
\sigma_k(u) x$ as $x\rightarrow\infty$ for some smooth function
$\sigma_k(u)$; then by Lemma 2.1
$$\frac{|P_{k+1}|}{\pi(x)} =\frac{\pi(x,P_k)}{\pi(x)}\sim \frac{\Psi(x,P_{k})}{x} \sim \sigma_k(u)
\hbox{ as } x\rightarrow \infty .$$ Now by corollary 1 we have
$$ \Psi(x,P_{k+1}) \sim \sigma_{k+1}(u) x \hbox{ as } x\rightarrow \infty, $$ where $\sigma_{k+1}(u)$
is the corresponding solution to (2) with $\chi(u) = \sigma_k(u)$.
Noting that $\sigma_0(u)= \rho(u) =((e+o(1))/u\log(u))^u$ and
using proposition 4 we deduce that
$$ \sigma_{k}(u)= \left(\frac{1+o(1)}{\log_{k}(u)\log_{k+1}(u)}\right)^u$$ by induction. Thus, using proposition 1,
the Theorem follows.
\enddemo

\head 3. Proof of Theorem 1 \endhead

\proclaim {Lemma 3.1} If $P$ is a set of primes $\leq x$, then
$$ \sum\Sb p\notin P\\ p\leq x \endSb\frac{1}{p} \ll 1 \qquad  \Longleftrightarrow  \qquad
\prod_{p\in P} \left(1-\frac{1}{p}\right) \asymp \frac{1}{\log(x)}.$$
\endproclaim
\demo{Proof}The result follows since
$$
\prod_{p\in
  P}\left(1-\frac{1}{p}\right)= \prod \Sb p\notin
  P\\ p\leq x \endSb \left(1-\frac{1}{p}\right)^{-1}\prod_{p\leq x}\left(1-\frac{1}{p}\right)\asymp
  \exp\left(O\left(\sum\Sb p\notin P\\ p\leq x \endSb\frac{1}{p}\right)\right)\frac{1}{\log(x)} $$by Mertens theorem.

\enddemo

\proclaim {Lemma 3.2} Let m,d be positive integers such that
$d|m$, then we have
$$ \sum\Sb r\leq x \\ d|r|m \endSb \frac{\mu(r)}{r} = \mu(d) \sum\Sb
n\geq 1 \\ d|n \\ p|n \implies p|d \endSb \frac{1}{n}\quad \sum\Sb r\leq
x/n \\ r|m \endSb \frac{\mu(r)}{r}.$$
\endproclaim

\demo{Proof}The result is trivial if $\mu(d)=0$ or $d=1$.
We fix $m$ and do a double induction on $d\geq 1$ and $x\geq 1$. Now

$$
\align
 S_d(x) &:= \sum\Sb r\leq x \\ d|r|m \endSb \frac{\mu(r)}{r} = \sum\Sb
n\leq x/d \\ n|\frac{m}{d} \endSb \frac{\mu(dn)}{dn} =  \frac{\mu(d)}{d}
\sum\Sb n\leq x/d\\n|m \\ (n,d)=1
\endSb\frac{\mu(n)}{n}=\frac{\mu(d)}{d} \sum\Sb n\leq x/d\\n|m
\endSb\frac{\mu(n)}{n}\sum\Sb a|n\\a|d \endSb\mu(a)\\ &=
\frac{\mu(d)}{d}\sum_{a|d}\mu(a)\sum\Sb n\leq x/d \\ a|n|m
\endSb \frac{\mu(n)}{n} = \frac{\mu(d)}{d}\sum_{a|d}\mu(a) S_a(x/d).
\endalign $$
Now each $a\leq d$ and $x/d < x$ so, by induction
$$
\align
S_d(x) &= \frac{\mu(d)}{d}\sum_{a|d}\mu(a)^2 \sum\Sb n\geq 1 \\ a|n \\
p|n \implies p|a \endSb \frac{1}{n} \sum\Sb r\leq x/nd \\ r|m
\endSb \frac{\mu(r)}{r}\\
& = \frac{\mu(d)}{d}\sum\Sb n\geq 1 \\ p|n
\implies p|d \endSb \frac{1}{n} \sum\Sb r\leq x/nd \\ r|m \endSb
\frac{\mu(r)}{r}\sum\Sb a|d\\a|n\\ p|n \implies p|a
\endSb \mu(a)^2.
\endalign
 $$
Now if we write $ n=p_1^{b_1}p_2^{b_2}....p_k^{b_k}$ with each $b_j \geq 1$, then $ p|n \implies p|d$ implies
that $ p_1p_2...p_k |d $. Moreover if  $a$ satisfies  $a|d$, $a|n$, $ p|n \implies p|a$, and $a$ is a
  squarefree, then $a$ must be $p_1p_2...p_k$; which implies $$ \sum\Sb a|d\\a|n\\ p|n \implies p|a
\endSb \mu(a)^2 =1 .$$
Then, writing $l=nd$, we have $\displaystyle{S_d(x) = \mu(d)\sum\Sb l\geq 1 \\ d|l \\ p|l
\implies p|d \endSb \frac{1}{l} \sum\Sb r\leq x/l \\ r|m
\endSb \frac{\mu(r)}{r}}$,  as desired.

\enddemo

\proclaim{Lemma 3.3} For any positive integer $k$ we have
$$ \sum\Sb n\geq 1\\ k|n\\ p|n \implies p|k\endSb \frac{\log(n)}{n} =
\frac{1}{\phi(k)}\left(\sum_{p|k}\frac{\log(p)}{p-1} + \log(k)\right)
\asymp \frac{\log(k)}{\phi(k)}.$$
\endproclaim

\demo{Proof} Writing $n= kd$ we have
$$ \quad \sum\Sb n\geq 1\\ k|n\\ p|n \implies p|k\endSb
\frac{\log(n)}{n} = \sum\Sb d\geq 1\\ p|d \implies p|k \endSb
\frac{\log(d)+\log(k)}{dk} = \frac{1}{k} \sum\Sb d\geq 1 \\ p|d
\implies p|k \endSb \frac{\log(d)}{d} + \frac{\log(k)}{\phi(k)}.
$$
Now if $ p_1,p_2,...,p_n$ are the prime factors of $k$ then
$$
\align
 &\sum\Sb d\geq 1\\ p|d \implies p|k \endSb \frac{\log(d)}{d} =
\sum\Sb a_i \geq 0 \\ 1\leq i\leq n\endSb
\frac{a_1\log(p_1)+a_2\log(p_2)+...+a_n\log(p_n)}{p_1^{a_1}p_2^{a_2}...p_n^{a_n}}\\
& =\sum_{i=1}^n\left(\sum_{a_i\geq 0}\frac{a_i\log(p_i)}{p_i^{a_i}}\right)\left(\prod\Sb 1\leq j \leq n
\\ j \not= i \endSb \left(\sum_{a_j\geq 0}
  \frac{1}{p_j^{a_j}}\right)\right)= \sum_{i=1}^n\frac{\log(p_i)}{p_i\left(1-\frac{1}{p_i}\right)^{2}}
  \prod\Sb 1\leq j \leq n \\ j \not= i \endSb\left(1-\frac{1}{p_j}\right)^{-1}\\ &= \frac{k}{\phi(k)}\sum_{p|k}
  \frac{\log(p)}{p-1}, \hbox{ which gives the result. }\\
\endalign
$$
\enddemo
We state a classical result of Sieve theory which is used
throughout the proof: \proclaim{Lemma 3.4 (Brun's Sieve)} Let $A$
be a set of positive integers contained in $[1,N]$. Suppose that
for each prime $p\leq N$, $A$ is excluded from $\omega(p)$ residue
classes$\mod p$, where $\omega$ is a multiplicative function and
$\omega(p)\ll 1$. Then
$$ |A| \ll N\prod_{p\leq N} \left(1-\frac{\omega(p)}{p}\right).$$
\endproclaim
\demo{Proof of Theorem 1}
Let $\epsilon >0$, $P^{*} = \{p\leq x\} \setminus P $, and $\displaystyle{ m = \prod\Sb
  p\in P^{*} \endSb p}$. Then we have
$$
\align
\pi(x,P)& = \sum\Sb p\leq x \\ q|p-1 \implies q \in P \endSb 1 =
\sum\Sb p\leq x \\ (p-1,m)=1 \endSb 1 = \sum_{p\leq x}\ \sum_{d|(m,p-1)}
\mu(d) = \sum_{d|m}\mu(d)\sum\Sb p\leq x\\ d|p-1\endSb 1\\ &= \sum_{d|m}\mu(d)\pi(x;d,1).\tag{3}\\
\endalign
$$Now by a similar argument we have
$$ \Psi(x,P) = \sum\Sb d\leq x\\  d|m\endSb
\mu(d)\left[\frac{x}{d}\right]. \tag{4}$$ By (3) and assuming
Conjecture 2 we have
$$\quad \pi(x,P) = \left(\sum\Sb d\leq x^{1-\epsilon} \\d|m\endSb
  \frac{\mu(d)}{\phi(d)}\right)\pi(x)+ O\left(\sum\Sb
  x^{1-\epsilon}<d\leq x\\ d|m \endSb \pi(x;d,1)\right) + o(\pi(x)) .\tag{5}$$
From (4), lemmas 3.1 and 3.4  we deduce
$$\left|\Psi(x,P)- x\sum\Sb d\leq x\\
d|m \endSb \frac{\mu(d)}{d}\right|\leq \sum\Sb d\leq x\\ d|m
\endSb 1 \leq \sum\Sb d\leq x \\ p|d \implies p\notin P \endSb
1\ll x\prod_{p\in P}\left(1-\frac{1}{p}\right)\ll
\frac{x}{\log(x)}.\tag{6}$$ Also by lemmas 3.1 and 3.4 we have
$$
\align
\sum\Sb x^{1-\epsilon}<d\leq x \\ d|m \endSb \frac{1}{d }& \leq\sum\Sb
x^{1-\epsilon}<d\leq x \\ p|d \implies p \in P^{*} \endSb \frac{1}{d} =
\int_{x^{1-\epsilon}}^x \frac{d\Psi(t,P^{*})}{t} \leq
\frac{\Psi(x,P^{*})}{x} + \int_{x^{1-\epsilon}}^x
\frac{\Psi(t,P^{*})}{t^2}dt\\
& \ll \prod_{p\in P}\left(1-\frac{1}{p}\right)\left(1+\int_{x^{1-\epsilon}}^x\frac{dt}{t}\right)\ll
\prod_{p\in P}\left(1-\frac{1}{p}\right) \epsilon \log(x)\ll
\epsilon.\tag{7}\\
\endalign
$$
Then from (5), (6) and (7) we deduce
$$
\align
\bigg|\frac{\pi(x,P)}{\pi(x)}&-\prod_{p\notin
  P}\left(1-\frac{1}{(p-1)^2}\right)\frac{\Psi(x,P)}{x}\bigg|\\ &  \leq
\left| \sum \Sb d\leq x^{1-\epsilon} \\d|m \endSb \frac{\mu(d)}{\phi(d)}-\sum \Sb d\leq x^{1-\epsilon} \\d|m
\endSb\frac{\mu(d)}{d}\prod_{p\notin
  P}\left(1-\frac{1}{(p-1)^2}\right)\right|\\ & + o(1) + O(\epsilon) + O\left(\sum\Sb x^{1-\epsilon}<d\leq x\\
  d|m \endSb
  \frac{\pi(x;d,1)}{\pi(x)}\right).\tag{8}
\endalign
$$
Now by lemmas 3.1, 3.4, and the fact that
$\displaystyle{\sum_{r\leq x}\frac{1}{\phi(r)}\ll \log(x)}$ we get
$$
\align
\quad \sum\Sb x^{1-\epsilon}<d\leq x\\ d|m \endSb
\pi(x;d,1) &= \sum_{r\leq x^{\epsilon}}\sum\Sb  x^{1-\epsilon} <d \leq x/r\\
p|d \implies p\notin P\\ dr+1 prime \endSb 1
 \ll \sum_{r\leq x^{\epsilon}}
\frac{x}{r}\prod_{p\in P} \left(1-\frac{1}{p}\right) \prod\Sb p\leq
x\\ p\nmid r \endSb \left(1-\frac{1}{p}\right) \\
& \ll \frac{x}{\log(x)^2}\sum_{r\leq x^{\epsilon}}\frac{1}{\phi(r)}
\ll \epsilon \frac{x}{\log(x)}. \tag{9}\\
\endalign
$$And from lemma 3.2 we have
$$
\align
 \sum\Sb d\leq x^{1-\epsilon}\\ d|m \endSb
\frac{\mu(d)}{\phi(d)} &= \sum\Sb d\leq x^{1-\epsilon}\\ d|m \endSb
\frac{\mu(d)}{d}\sum_{k|d}\frac{\mu(k)^2}{\phi(k)}\\ & =
\sum_{k|m}\frac{\mu(k)^2}{\phi(k)}\sum\Sb d\leq x^{1-\epsilon}\\ k|d|m
\endSb \frac{\mu(d)}{d} = \sum_{k|m}\frac{\mu(k)}{\phi(k)} \sum\Sb
 n\geq 1\\ k|n \\ p|n \implies p|k \endSb \frac{1}{n} \sum\Sb r
\leq  x^{1-\epsilon}/n \\ r|m \endSb \frac{\mu(r)}{r}\\
& =  \sum_{k|m}\frac{\mu(k)}{\phi(k)} \sum\Sb
 n\geq 1\\ k|n \\ p|n \implies p|k  \endSb \frac{1}{n} \sum\Sb r
\leq  x^{1-\epsilon} \\ r|m \endSb \frac{\mu(r)}{r} - \sum_{k|m}\frac{\mu(k)}{\phi(k)} \sum\Sb
 n\geq 1\\ k|n \\ p|n \implies p|k  \endSb \frac{1}{n} \sum\Sb x^{1-\epsilon}/n<r
\leq  x^{1-\epsilon} \\ r|m \endSb \frac{\mu(r)}{r}. \tag{10}\\
\endalign
$$
The first term in the RHS of (10) is equal to:
$$\sum\Sb r\leq x^{1-\epsilon}\\r|m \endSb
\frac{\mu(r)}{r}\sum_{k|m} \frac{\mu(k)}{k\phi(k)}
\prod_{p|k}\left(1-\frac{1}{p}\right)^{-1} = \sum\Sb r\leq
x^{1-\epsilon} \\ r|m \endSb
\frac{\mu(r)}{r}\prod_{p|m}\left(1-\frac{1}{(p-1)^2}\right).\tag{11}$$
By integration by parts and using lemma 3.4 we have
$$
\align
\left|\sum\Sb x^{1-\epsilon}/n <r \leq x^{1-\epsilon}\\ r|m
\endSb \frac{\mu(r)}{r}\right|& \leq \sum\Sb x^{1-\epsilon}/n <r \leq
x^{1-\epsilon}\\ r|m \endSb \frac{1}{r} \leq
\int_{x^{1-\epsilon}/n}^{x^{1-\epsilon}}\frac{d\Psi(t,P^{*})}{t}\\& \ll
\prod_{p\in P}\left(1-\frac{1}{p}\right)\left(1+\int_{x^{1-\epsilon}/n}^{x^{1-\epsilon}}\frac{dt}{t}\right)
\ll \frac{\log(n)}{\log(x)}.\\
\endalign
$$ Then, by lemma 3.3
$$
\align
 \sum_{k|m}\frac{\mu(k)}{\phi(k)} \sum\Sb
 n\geq 1\\ k|n \\ p|n \implies p|k \endSb \frac{1}{n} \sum\Sb x^{1-\epsilon}/n<r
\leq  x^{1-\epsilon} \\ r|m \endSb \frac{\mu(r)}{r}
&\ll \sum_{k|m}\frac{\mu(k)}{\phi(k)}
  \sum \Sb n\geq 1 \\ k|n \\ p|n \implies p|k \endSb
  \frac{\log(n)}{n\log(x)}\\& \ll \frac{1}{\log(x)} \sum_{k|m}
  \frac{\mu(k)\log(k)}{\phi(k)^2} \ll \frac{1}{\log(x)}. \tag{12}\\
\endalign $$
Thus combining (8), (9), (10), (11) and (12) gives the result,
letting $\epsilon \to 0$.

\enddemo

\head 4. Proof of Proposition 1 \endhead
\proclaim {Lemma 4.1}
$P_0=\{p\leq y\} \subseteq P_1 \subseteq P_2\subseteq...\subseteq
P_k \subseteq...$ where $P_{k+1} = \{$primes $q\leq x : p|q-1
\implies p \in P_k\}$.
\endproclaim
\demo {Proof} If $p \in P_0$ then $p \leq y$ and so $p-1 \leq y$,
which implies $ q|p-1 \implies q \leq y$. This means that $p \in
P_1$. Now using a simple induction argument: if $p \in P_k$ then
$q|p-1 \implies q \in P_{k-1} \subseteq P_k $, and so $p\in
P_{k+1}$.
\enddemo
\proclaim{Lemma 4.2}Let $r$ be a positive integer. Then
$$ R(r,k,x) := \sum\Sb y<r<q_1<...<q_k\leq x\\r|q_1-1,q_1|q_2-1,...,q_{k-1}|q_k-1\endSb \frac{1}{q_k} \leq
\frac{(\log x+1)^k}{r}.$$
We deduce that
$$ S(r,k,x) := \sum\Sb y<r<q_1<...<q_k\leq n\leq  x\\r|q_1-1,q_1|q_2-1,...,q_{k-1}|q_k-1,q_k|n\endSb 1 \leq
\frac{x(\log x+1)^k}{r}.$$
\endproclaim
\demo{Proof}
Writing $q_k -1 = mq_{k-1}$ we have
$$ R(r,k,x)\leq \sum_{m\leq \frac{x}{r}}\frac{1}{m}R(r,k-1,x)\leq R(r,k-1,x)(\log x+1),$$and $$R(r,1,x)\leq
\sum_{m\leq \frac{x}{r}} \frac{1}{mr} \leq \frac{\log x+1}{r}, $$ then by induction $$ R(r,k,x) \leq
\frac{(\log x+1)^k}{r}.$$ The second inequality follows since
$$S(r,k,x)\leq  \sum\Sb y<r<q_1<...<q_k\leq x\\r|q_1-1,q_1|q_2-1,...,q_{k-1}|q_k-1\endSb \frac{x}{q_k}=x R(r,k,x).$$

\enddemo
\proclaim {Lemma 4.3} Define $$S_k(x,y) = \{n\leq x : \hbox{ there
is a prime } p > y \hbox{ such that } p^2|\phi_k(n) \}$$ Then
$$\left|\Psi(x,P_k)-\Phi_k(x,y)\right| \leq \sum_{i=0}^{k-1}|S_i(x,y)|.$$
\endproclaim
\demo{Proof} Let $A_k(x) = \{n\leq x : p|n \implies p\in P_k\}$.
If $n\in A_{k+1}(x)$ and $\phi(n)\notin A_k(x)$, then there is a
prime $p$ which divides $\phi(n)$ and $p\notin P_k$. Now $n\in
A_{k+1}(x)$ so every prime factor of $q-1$, where $q|n$, is in
$P_k$, which implies that $p^2|n$. This gives
$$ A_{k+1}(x) \setminus \{n\leq x : \phi(n) \in A_k(x)\} = \{n\leq x: n\in A_{k+1}(x) ,
\exists\hbox{ a prime }p\in P_{k+1}\setminus P_k, p^2|n\}. $$
Then by lemma 4.1
$$
\align
0&\leq\Psi(x,P_k)-\Phi_k(x,y)=|A_k(x)|-|\{n\leq x : \phi_k(n) \in A_0(x)\}|\\
&=\sum_{i=0}^{k-1}\left|\{n\leq x : \phi_i(n) \in A_{k-i}(x)\}\right|-\left|\{n\leq x : \phi_{i+1}(n)
\in A_{k-i-1}(x)\}\right|
\\ &= \sum_{i=0}^{k-1} |\{n\leq x: \phi_{i}(n)\in A_{k-i}(x) ,
\hbox{ there is a prime } p\in P_{k-i}\setminus P_{k-i-1}, p^2|\phi_{i}(n) \}|\\&\leq \sum_{i=0}^{k-1}|S_i(x,y)|.\\
\endalign
$$
\enddemo
\demo{Proof of Proposition 1}
Note that if $q|(\phi(n),n)$ for some prime $q$, then $q^2|n$. Define
$$S_k^*(x,y) = S_k(x,y) \setminus \bigcup_{i=0}^{k-1}S_i(x,y).$$
If $n \in S_k^*(x,y)$ and $q^2|\phi_j(n)$ for some $0\leq j\leq
k-1$, then $q\leq y$ (by definition); also there exists some prime
$p$ satisfying $p^2|\phi_k(n)$ with $p > y$, which implies $p^2
\nmid \phi_{k-1}(n)$. Thus we have two cases :

(i) There exists a prime $q_1|\phi_{k-1}(n)$ such that $p^2|q_1-1$ .

(ii) There are two primes $q_1|\phi_{k-1}(n)$ and $Q_1|\phi_{k-1}(n)$ such that $p|q_1-1$ and $p|Q_1-1$.

In the first case $q_1|\phi_{k-1}(n)=\phi(\phi_{k-2}(n))$,
$p|q_1-1$, so that $q_1 > y$, which implies that $q_1^2\nmid
\phi_{k-2}(n)$, so that there exists a prime $q_2|\phi_{k-2}(n)$
such that $q_1|q_2-1$ and $q_2>q_1>p> y$. By a simple induction,
there exist primes $y <p<q_1<q_2<...<q_k$ for which $
p^2|q_1-1,q_1|q_2-1,...,q_{k-1}|q_k-1,q_k|n$.

\noindent We deduce that the total number of possibilities for this case is:
$$
S_1=\sum\Sb y<p<q_1<...<q_k\leq n\leq x\\p^2|q_1-1,q_1|q_2-1,...,q_{k-1}|q_k-1,q_k|n\endSb 1 =
 \sum_{y<p<\sqrt{x}}S(p^2,k,x) \leq  x(\log x+1)^k\sum_{p>y}\frac{1}{p^2} \ll \frac{x(\log x)^k}{y}
$$
by lemma 4.2.

Now following an analogous argument we find (for the second case)
that there exist primes $p,q_1,q_2,...,q_k,Q_1,Q_2,...,Q_k$ such
that $ p|q_1-1,q_1|q_2-1,...,q_{k-1}|q_k-1,q_k|n$ and $
p|Q_1-1,Q_1|Q_2-1,...,Q_{k-1}|Q_k-1,Q_k|n$; we'll have two cases
again:

a) $q_i \neq Q_i$ for all $1 \leq i \leq k$.

b) There exists $i$ such that  $q_i= Q_i$; so let  $j= \min\{1\leq i\leq k :q_i= Q_i\}$.

\noindent For case a) the total number of possibilities is:
$$
\align
S_2&=\sum\Sb y<p<q_1<...<q_k\leq n\leq x\\p<Q_1<...<Q_k\leq n \leq x
\\p|q_1-1,q_1|q_2-1,...,q_{k-1}|q_k-1,q_k|n\\p|Q_1-1,Q_1|Q_2-1,...,Q_{k-1}|Q_k-1,Q_k|n\endSb 1
\leq \sum\Sb y<p<q_1<...<q_k\leq x\\p<Q_1<...<Q_k\leq x\\p|q_1-1,q_1|q_2-1,...,q_{k-1}|q_k-1\\p|Q_1-1,Q_1|Q_2-1,...,
Q_{k-1}|Q_k-1\endSb \frac{x}{q_kQ_k}\\
& \leq x \sum_{p> y} R(p,k,x)^2
=O\left(\frac{x(\log x)^{2k}}{y}\right) \\
\endalign
$$
by lemma 4.2.

\noindent Now for case b) $ p|q_1-1,q_1|q_2-1,...,q_{j-1}|q_j-1$,
$ p|Q_1-1,Q_1|Q_2-1,...,Q_{j-1}|Q_j-1$ and
$Q_j=q_j|\phi_{k-j}(n)$, then following the same logic there exist
primes $q_{j+1},q_{j+2},...,q_k$ such that $q_j|q_{j+1}-1
,...,q_{k-1}|q_k-1,q_k|n$.

\noindent We deduce that the total number of possibilities is:
$$
\align
S_3&=\sum\Sb y<p<q_1<...<q_j\leq x\\p|q_1-1,q_1|q_2-1,...,q_{j-1}|q_j-1\endSb \left(\sum\Sb p<Q_1<...<Q_j=q_j
\\ p|Q_1-1,Q_1|Q_2-1,...,Q_{j-1}|Q_j-1\endSb\left(\sum \Sb qj<q_{j+1}<...<q_k\leq n \leq x
\\ q_j|q_{j+1}-1,...q_{k-1}|q_k-1,q_k|n\endSb 1\right)\right)\\
&=\sum\Sb y<p<q_1<...<q_j\leq x\\p|q_1-1,q_1|q_2-1,...,q_{j-1}|q_j-1\endSb\left(\sum \Sb p<Q_1<...<Q_j=q_j
\\p|Q_1-1,Q_1|Q_2-1,...,Q_{j-1}|Q_j-1\endSb S(q_j,k-j,x) \right)\\
&\leq \sum\Sb y<p<q_1<...<q_j\leq x\\p|q_1-1,q_1|q_2-1,...,q_{j-1}|q_j-1\endSb\left(\sum \Sb p<Q_1<...<Q_j=q_j
\\p|Q_1-1,Q_1|Q_2-1,...,Q_{j-1}|Q_j-1\endSb \frac{x(\log x+1)^{k-j}}{q_j}\right).\\
\endalign
$$
Now writing $Q_j-1=q_j-1=mQ_{j-1}q_{j-1}$ we have:
$$
\sum\Sb q_j\leq x\\q_{j-1},Q_{j-1}|q_j-1\endSb \frac{1}{q_j} <
\sum_{m\leq x}\frac{1}{mQ_{j-1}q_{j-1}} \leq \frac{\log
x+1}{Q_{j-1}q_{j-1}}.$$ Thus by lemma 4.2
$$
S_3 \ll x(\log x)^{k-j+1}\sum_{p>y}R(p,j-1,x)^2\ll \frac{x(\log x)^{k+j-1}}{y}.
$$
We deduce from cases (i), (ii) a) and b) that
$$ |S_k^*(x,y)| = S_1+ S_2 +S_3 =  O\left(\frac{x(\log x)^{2k}}{y}\right).\tag {13}$$
Now
$$
\align
|S_1(x,y)|&= |S_0(x,y)|+|S_1^*(x,y)|=|\{n\leq x :\exists \hbox{  a prime } p > y, p^2|n\}| + O\left(\frac{x(\log x)^2}{y}
\right)\\
& \leq \sum\Sb n\leq x\\ p> y\\p^2|n \endSb 1 +O\left(\frac{x(\log x)^2}{y}\right) \leq \sum_{p> y}\frac{x}{p^2}
+O\left(\frac{x(\log x)^2}{y}\right)\\
& = O\left(\frac{x}{y}\right)+O\left(\frac{x(\log x)^2}{y}\right) = O\left(\frac{x(\log x)^2}{y}\right),\\
\endalign
$$
and by simple induction we obtain:
$$|S_k(x,y)| = |S_k^*(x,y)| + \sum_{i=0}^{k-1}|S_i(x,y)| = O\left(\frac{x(\log x)^{2k}}{y}\right). \tag{14}$$
Thus by (14) and lemma 4.3 the result follows.

\enddemo

\head 5. Proof of Theorem 3\endhead

\proclaim{Lemma 5.1} Let $\chi$ be a real measurable function for
which $\int_1^{\infty}\chi(t) e^{\xi t} dt$ converges for all
$\xi$ and such that $C:=\int_1^{\infty}\chi(v) dv >0$. Then for
$u\geq C^2$ and for any $\epsilon >0$, we have
$$\int_1^{\infty}\chi(v) e^{(\xi(u)+\epsilon) v}
dv\geq u^{1+\frac{\epsilon}{2\xi(u)}}.$$
\endproclaim

\demo{Proof} Let $\epsilon'>0$ and $s>0$. Using H\"older
inequality we get
$$ \left(\int_1^{\infty}\chi(v)dv\right)^{\epsilon'}\left(\int_1^{\infty}\chi(v) e^{s v}
dv\right)^{1-\epsilon'}\geq \int_1^{\infty}\chi(v)
e^{s(1-\epsilon') v} dv.$$ Now putting
$s=\xi(u)/(1-\epsilon')=\xi(u)(1+\epsilon")$ and since $u\geq C^2$
we deduce that
$$\int_1^{\infty}\chi(v) e^{\xi(u)(1+\epsilon") v} dv\geq \frac{1}{C^{\epsilon'}}
\left(\int_1^{\infty}\chi(v) e^{\xi(u) v}
dv\right)^{1+\epsilon"}\geq u^{1+\epsilon"/2}.$$ The lemma follows
taking $\epsilon=\xi(u)\epsilon"$.
\enddemo

\demo{Proof of Theorem 3} \enddemo \demo{The Upper bound}

 From Lemma 3.4 of Granville-Soundararajan [8]
we note that
$$
\sigma(u) = \rho(u) + \sum_{j=1}^{\infty} \frac{1}{j!} \int\Sb
t_1, \ldots, t_j \ge 1\\ t_1+\cdots + t_j \le u\endSb
\frac{\chi(t_1)}{t_1} \ldots \frac{\chi(t_j)}{t_j} \rho(u-t_1
-\ldots-t_j) dt_1\ldots dt_j.
$$
Therefore for any $\xi \in {\Bbb R}$
$$
\align \sigma(u) e^{\xi u} &= \rho(u) e^{\xi u} +
\sum_{j=1}^{\infty} \frac{1}{j!}
\\
&\int\Sb t_1, \ldots, t_j \ge 1\\ t_1+\cdots + t_j \le u\endSb
\frac{\chi(t_1)e^{\xi t_1}}{t_1} \ldots \frac{\chi(t_j)e^{\xi
t_j}}{t_j} \rho(u-t_1 -\ldots-t_j) e^{\xi(u-t_1-\ldots -t_j)}
dt_1\ldots dt_j.\\
\endalign
$$
Setting $F(\xi)=\max_{t \ge 0} \rho(t) e^{\xi t}$ we deduce that
(by forgetting the condition $t_1+\ldots+t_j \le u$)
$$
\sigma(u) \le F(\xi) e^{-\xi u} \sum_{j=0}^{\infty} \frac{1}{j!}
\Big(\int_1^{\infty} \frac{\chi(t) e^{\xi t}}{t} dt \Big)^j =
F(\xi) e^{-\xi u} \exp\Big(\int_1^{\infty} \frac{\chi(t) e^{\xi
t}}{t} dt\Big).
$$
Choose $\xi$ such that $u=\int_1^{\infty} \chi(t) e^{\xi t} dt$,
that is $\xi=\xi(u)$.

\noindent Now putting $C:=\int_1^{\infty}\chi(v) dv$ we have   $
u=\int_1^{\infty}\chi(v) e^{\xi(u) v} dv\geq  C e^{\xi(u)}$, which
implies that
$$ F(\xi(u))\leq \max_{t \ge 0} \left(\frac{(e+o(1))u}{t\log(t)C}\right)^t=e^{O\left(u/\log(u)\right)},$$
and the upper bound follows.
\enddemo
\demo{The Lower bound }Fix $\epsilon>0.$

We will show that there exists a constant $C_{\epsilon}$ such that

$$\sigma(u)>C_{\epsilon}\exp\left((-\xi(u)-\epsilon) u+\int_1^{\infty}\frac{\chi(v) e^{\xi(u) v}}{v}
dv\right)\hbox { for all } u\geq 0.\tag{15}$$ Let $u_0$ be a
suitably large number, and define
$$C_{\epsilon}=C_{\epsilon,u_0}=\inf_{u\leq u_0}\sigma(u)\exp\left((\xi(u)+\epsilon) u-\int_1^{\infty}\frac{\chi(v)
 e^{\xi(u) v}}{v}
dv\right).$$ Evidently (15) holds for all $u\leq u_0$.

We use an induction argument. Let $n\in \Bbb N$ such that $n>u_0$
and suppose that (15) is verified for all $t\leq n$, then we will
show that (15) holds for all $ t\in [n,n+1]$.

Define $f(\xi)=\displaystyle{\int_1^{\infty}\frac{\chi(v) e^{\xi
v}}{v} dv}$, and let $u\in [n,n+1]$. Then using our hypothesis we
have
$$\align
&\frac{\displaystyle{\sigma(u)e^{\left((\xi(u)+\epsilon) u\right)}}}{\displaystyle{C_{\epsilon}\exp(f(\xi(u)))}}=
\frac{1}{C_{\epsilon}u}\displaystyle{\int_0^u\chi(t)e^{(\xi(u)+\epsilon)t}\sigma(u-t)
e^{((\xi(u)+\epsilon)(u-t)-f(\xi(u))}dt}\\
&\geq
\displaystyle{\frac{1}{u}\int_1^u\chi(t)\exp\Big((\xi(u)+\epsilon)t+
(\xi(u)-\xi(u-t))(u-t)+f(\xi(u-t))-f(\xi(u))\Big)dt.}\tag{16}
\endalign$$
Since $f'(\xi)=\int_1^{\infty}\chi(v) e^{\xi v} dv$ and using the
mean value theorem we deduce that
$$ u-t \leq \frac{f(\xi(u))-f(\xi(u-t))}{\xi(u)-\xi(u-t)}\leq u.\tag{17}$$
Now differentiating $u=\int_1^{\infty}\chi(v) e^{\xi(u) v} dv$
with respect to $u$ we get that
$$\xi'(u) = \displaystyle{\left(\int_1^{\infty}v\chi(v)e^{\xi(u) v}dv\right)^{-1}}\leq \frac{1}{u}.\tag{18}$$
By (18) and using the mean value theorem again we have
$$\xi(u)-\xi(u-t)\leq \frac{t}{u-t}.\tag{19}$$
Therefore by (17), then (19) we deduce that
$$
\align
\frac{1}{u}\int_1^u\chi(t)&\exp\Big((\xi(u)+\epsilon)t+(\xi(u)-\xi(u-t))(u-t)+f(\xi(u-t))-f(\xi(u))\Big)dt\\
&\geq \displaystyle{\frac{1}{u}\int_1^u\chi(t)\exp\Big((\xi(u)+\epsilon)t-t(\xi(u)-\xi(u-t))\Big)dt}\\& \geq
\displaystyle{\frac{1}{u}\int_1^u\chi(t)\exp\left((\xi(u)+\epsilon)t-\frac{t^2}{(u-t)}\right)dt}\\&  \geq
\displaystyle{\frac{1}{u}\int_1^{\sqrt{u}}\chi(t)e^{(\xi(u)+\epsilon/2)t}dt,} \ \ \text{\rm for} \ \ u \geq u_0 .
\tag{20}\\
\endalign$$
For case i), Since $\int_T^{\infty}\chi(t)dt=0$ and $\chi(t)\geq
0$ for all $t$, then meas$\{t\geq T: \chi(t)\ne 0 \}=0$ which
implies that meas$\{t\geq T: \chi(t)e^{\xi(u)t}\ne 0 \}=0$, and so
$\int_T^{\infty}\chi(t)e^{\xi(u)t}dt=0$. Then taking $u_0>T^2$ we
have
$$ \int_1^{\sqrt{u}}\chi(t)e^{(\xi(u)+\epsilon/2)t}dt = \int_1^{T}\chi(t)e^{(\xi(u)+\epsilon/2)t}dt > \int_1^{T}
\chi(t)e^{\xi(u)t}dt = u . \tag{21}$$ Now for case ii) since
$\chi(t)=\left(\{1+o(1)\}/h(t)\right)^t$, there exist two
constants $A_\epsilon$ and $B_\epsilon$ for which
$$ A_\epsilon \left(\frac{\exp(-\epsilon/16)}{h(t)}\right)^t  <\chi(t) < B_\epsilon \left(\frac{\exp(\epsilon/16)}
{h(t)}\right)^t \ \ \text{\rm for every} \ \ t\geq 0.\tag{22}$$ we
consider two cases

\enddemo
\demo{1)} \qquad $\displaystyle{\frac{e^{\xi(u)}}{h(\sqrt{u})}\geq
\exp\left(-\frac{\epsilon}{4}\right).}$ Since $h$ is
non-decreasing we have by (22)
$$
\align
\int_1^{\sqrt{u}}\chi(t)e^{(\xi(u)+\epsilon/2)t}dt &\geq A_\epsilon \int_1^{\sqrt{u}}
\left(\frac{e^{(\xi(u)+7\epsilon/16)}}{h(t)}\right)^tdt \geq A_\epsilon\int_1^{\sqrt{u}}
\left(\frac{e^{\xi(u)}}{h(\sqrt{u})}\right)^t e^{(7\epsilon/16)t}dt\\
&\geq A_\epsilon\int_1^{\sqrt{u}}e^{(3\epsilon/16)t}dt = \frac{16A_{\epsilon}}{3\epsilon}
\left(e^{(3\epsilon/16)\sqrt{u}}-e^{3\epsilon/16}\right)>u ,\tag{23}\\
\endalign
$$for $u>u_0$.
\enddemo

\demo{2)} \qquad $\displaystyle{\frac{e^{\xi(u)}}{h(\sqrt{u})}\leq
\exp\left(-\frac{\epsilon}{4}\right).}$ Using lemma 5.1 and (22),
then the fact that $h$ is non-decreasing and $u \geq Ce^{\xi(u)}$
we conclude that
$$
\align
\int_1^{\sqrt{u}}\chi(t)e^{(\xi(u)+\epsilon/2)t}dt &\geq \int_1^{\sqrt{u}}\chi(t)e^{(\xi(u)+\epsilon/8)t}dt
\\ &\geq u^{1+\epsilon/(16\xi(u))} - B_\epsilon\int_{\sqrt{u}}^{\infty}
\left(\frac{e^{(\xi(u)+3\epsilon/16)}}{h(t)}\right)^tdt\\
&\geq u^{1+\epsilon/(16\xi(u))} - B_\epsilon\int_{\sqrt{u}}^{\infty}\left(\frac{e^{\xi(u)}}{h(\sqrt{u})}\right)^t
 e^{(3\epsilon/16)t}dt\\ &\geq u^{1+\epsilon/(16\xi(u))} -
 B_\epsilon\frac{16}{\epsilon}\exp\Big(-\epsilon\sqrt{u}/16\Big) >u, \tag{24} \\
\endalign
$$for $u>u_0$.
Thus using (16), (20) then (21), (23), and (24) the result
follows.
\enddemo

\head 6. Getting the asymptotic of $\sigma$ explicitly \endhead

\proclaim{Lemma 6.1} If $\xi(u)=o(\log(u))$ as $u\rightarrow
\infty$, then
$$\int_1^{\infty}\frac{\chi(v)e^{\xi(u)v}}{v}dv = o(u) \ \
\text{\rm and so} \ \ \sigma(u)=\exp((-\xi(u)+o(1))u).$$
\endproclaim
\demo{Proof} Since $\chi(t)\leq 1$ for every $t\geq 1$ and using
our assumption we have
$$
\align
\int_1^{\infty}\frac{\chi(v)e^{\xi(u)v}}{v}dv &= \int_1^{\frac{\log(u)}{\xi(u)}}\frac{\chi(v)e^{\xi(u)v}}{v}dv +
\int_{\frac{\log(u)}{\xi(u)}}^{\infty}\frac{\chi(v)e^{\xi(u)v}}{v}dv\\
& \leq \int_1^{\frac{\log(u)}{\xi(u)}}e^{\xi(u)v}dv +
\frac{\xi(u)}{\log(u)}\int_1^{\infty}\chi(v)e^{\xi(u)v}dv\\
&=
\frac{1}{\xi(u)}\left(u-e^{\xi}\right)+\frac{\xi(u)u}{\log(u)}=o(u).
\endalign
$$
\enddemo

\demo{Proof of Proposition 3} Let $\zeta(u)$ be the unique
continuous solution to the equation $u=e^{\zeta(u)T}/\zeta(u)$.
Since $\chi(t)\leq 1$ for all $t$, we have

$$\frac{e^{\zeta(u)T}}{\zeta(u)}= \int_1^{T}\chi(v)e^{\xi(u)v}dv\leq \frac{e^{\xi(u)T}-e^{\xi(u)}}{\xi(u)}<
\frac{e^{\xi(u)T}}{\xi(u)},$$ and since the function
$f(\xi)=e^{\xi T}/\xi$ is non-decreasing for $\xi>1$ we deduce
that $\zeta(u)\leq\xi(u)$. Now fix $\epsilon>0$ (such that
$T(1-\epsilon)>1$), and suppose that there is arbitrary large $u$
for which $\xi(u)>\zeta(u)(1+\epsilon)$. Define $s_\epsilon=\int_{
T(1-\epsilon/3)}^T\chi(t)dt>0$ (by the definition of $T$) . We
deduce under our assumption that
$$s_\epsilon e^{\zeta(u)(1+\epsilon/3)T}\leq  s_\epsilon e^{T(1-\epsilon/3)\xi(u)}
\leq\int_{T(1-\epsilon/3)}^{T}\chi(v)e^{\xi(u)v}dv \leq
\frac{e^{\zeta(u)T}}{\zeta(u)},
$$which is impossible if $u$ is large enough.
 Thus
$ \xi(u) = \zeta(u)(1+o(1))$ as $u\rightarrow \infty$. Now we
trivially have $1\ll \zeta(u) \ll \log(u)$, then
$$\align
\zeta(u) &= \frac{\log(u)}{T}+\frac{\log(\zeta(u))}{T}=\frac{\log(u)}{T}(1+o(1)).\\
\endalign$$ We deduce that $\xi(u)=\displaystyle{\frac{\log(u)}{T}}(1+o(1))$,
 and the result follows combining Theorem 3 and the fact that $\displaystyle{\int_1^{T}
 \frac{\chi(v)e^{\xi(u)v}}{v}dv}=O(u)$.

\enddemo
\noindent Now we prove Proposition 4; define
$g(u):=h(u)/(uh'(u))$.

 \proclaim{Lemma 6.2} Let $h(u)$ be a real differentiable function
with $uh'(u)/h(u)= n+o(1)$, where $n$ is a positive constant. Then
for all $ k>0$ we have $h(ku)=h(u) k^{n +o(1)}$.
\endproclaim

\demo {Proof} We have that
$$\log \left(\frac{h(ku)}{h(u)}\right) = \int_u^{ku}\frac{h'(t)}{h(t)}dt
=\int_u^{ku}\frac{(n+o(1))}{t}dt = (n+o(1)) \log k.$$
\enddemo

\proclaim{Lemma 6.3} Assume the hypothesis of Proposition 4 b).

 \noindent Then
$h(v(u)\log(u))= (1+o(1))h(\log u)$, where $v(u) :=
\displaystyle{\min\left(\log(u),\min_{\log(u)\leq t \leq
\log^2(u)}g(t)\right)}$, and $v(u)\rightarrow \infty$ as
$u\rightarrow \infty$.
\endproclaim

\demo{Proof} Since $g(t)\rightarrow \infty$ as $t\rightarrow
\infty$ then $v(u)\rightarrow \infty$ as $u\rightarrow \infty$, so
if $u$ is large then $v(u)\log(u)>\log u$ so that
$h(v(u)\log(u))\geq h(\log u)$. On the other hand
$$
\align &\log\left(\frac{h(v(u)\log(u))}{h(\log(u))}\right)
 = \int_{\log(u)}^{v(u)\log(u)}\left(\frac{h'(t)}{h(t)}\right) dt
= \int_{\log(u)}^{v(u)\log(u)}\frac{1}{tg(t)}dt\\
& \leq \frac{1}{\min_{\log(u)\leq t \leq v(u)\log(u)}g(t)}
\int_{\log(u)}^{v(u)\log(u)}\frac{dt}{t}
\leq  \frac{1}{v(u)}(1+ \log(v(u))) = o(1) .\\
\endalign
$$
\enddemo

\demo{Proof of Proposition 4} Fix $\epsilon>0$ and suppose that
there is arbitrary large $u$ for which $\xi(u) >
\log\left(h(\zeta\log(u))\right)+\epsilon$. Then for such $u$ we
have:

\noindent In case a) by (22), lemma 6.2 and since $h$ is
non-decreasing
$$
\align
 u &=\int_1^{\infty}\chi(t) e^{\xi(u) t}
dt>
A_\epsilon\int_1^{\infty}\left(\frac{e^{\xi(u)-\epsilon/16}}{h(t)}\right)^tdt
\geq A_\epsilon\int_{\log\log(u)}^{\log(u)/n}e^{(\epsilon/2)
t}\left(\frac{h(e\log(u)/n)}{h(\log(u)/n)}\right)^tdt\\ & =
A_\epsilon\int_{\log\log(u)}^{\log(u)/n}e^{(\epsilon/2
+n+o(1))t}dt >
A_\epsilon\int_{\log\log(u)}^{\log(u)/n}e^{(\epsilon/4 +n)t}dt>u,
\endalign
$$for $u$ large enough, which is a contradiction.

\noindent Now in case b), our assumption and lemma 6.3 imply that
$\xi(u)>\log\left(h(v(u)\log(u))\right)+\epsilon/2$. Then by (22)
and since $h$ is non-decreasing and $v(u)\rightarrow \infty$ as
$u\rightarrow \infty$

$$
\align u &= \int_1^{\infty}\chi(t) e^{\xi(u) t}
dt> A_\epsilon\int_1^{\infty}\left(\frac{e^{\xi(u)-\epsilon/16}}{h(t)}\right)^tdt \geq A_\epsilon
\int_{1}^{v(u)\log(u)}e^{(\epsilon/3) t}dt\\
 &= A_\epsilon\frac{3}{\epsilon}\left(u^{v(u)\epsilon/3}-e^{\epsilon/3}\right)>u,\endalign$$
 for $u$ large enough, which is a contradiction.

Now we suppose that there is arbitrary large $u$ for which $\xi(u)
< \log\left(h(\zeta\log(u))\right)-\epsilon$. Then for such $u$
let $q(t):=(\xi(u)+\epsilon/16-\log(h(t)))t$, so that $$
q'(t)=\xi(u)+
\frac{\epsilon}{16}-\log(h(t))-\displaystyle{\left(\frac{th'(t)}{h(t)}\right)}.$$
Now in case a) $ q'(t)=\xi(u)+\epsilon/16-\log(h(t))- n+o(1)$,
therefore the maximum of $q(t)$ holds at some point $t_0$ for
which $q'(t_0)=0$ so that, under our assumption
$$h(t_0)=e^{\xi(u)+\epsilon/16-n+o(1)} < h(e\log(u)/n)
e^{-n-\epsilon/2}.\tag{25}$$ Now we must have
$$t_0<\log(u)(1-\epsilon/(4n))/n,\tag{26}$$ otherwise since $h$ is
non-decreasing and  by lemma 6.2
$$ \frac{h(t_0)}{h(e\log(u)/n)} \geq \frac{h\left(\frac{\log(u)}{n}\left(1-\frac{\epsilon}{4n}\right)\right)}
{h(e\log(u)/n)}=\left(\left(1-\frac{\epsilon}{4n}\right)e^{-1}\right)^{n+o(1)}>e^{-n-\epsilon/2}
\ \ \text{\rm contradicting (25)}.$$ By (22), (25) and (26) and
since $h$ is non-decreasing we deduce that
$$
\align u &=\int_1^{\infty}\chi(t) e^{\xi(u) t}
dt < B_\epsilon\int_1^{\infty}\left(\frac{e^{\xi(u)+\epsilon/16}}{h(t)}\right)^tdt\\ & =
B_\epsilon\int_1^{e\log(u)/n}\left(\frac{e^{\xi(u)+\epsilon/16}}{h(t)}\right)^tdt+
B_\epsilon\int_{e\log(u)/n}^{\infty}\left(\frac{e^{\xi(u)+\epsilon/16}}{h(t)}\right)^tdt\\
& \leq B_\epsilon\frac{e\log(u)}{n}
\left(\frac{e^{\xi(u)+\epsilon/16}}{h(t_0)}\right)^{t_0} +
B_\epsilon\int_{e\log(u)/n}^{\infty}e^{-\epsilon/2 t} dt =
B_\epsilon\frac{e\log(u)}{n}(e^{n+o(1)})^{t_0} +o(1) <u,
\endalign$$
for $u$ large enough, which is a contradiction.

\noindent For case b) $
q'(t)=\xi(u)+\epsilon/16-\log(h(t))-\frac{1}{g(t)},$ and the
maximum of $q(t)$ holds at some point $t_0$ for which $q'(t_0)=0$
(to avoid redundancy we take $t_0=\min\{t: q'(t)=0\}$, which is
possible by the continuity of $h(t)$ and $g(t)$). Now
$t_0\rightarrow \infty$ as $u\rightarrow \infty$, otherwise
$q'(t_0)>0$ for $u$ large enough. Thus
$$\left(\frac{e^{\xi(u)+\epsilon/16}}{h(t_0)}\right)^{t_0}=\exp\left(\frac{t_0}{g(t_0)}\right)=e^{o(t_0)}.\tag{27}$$
Now by (22)
$$
\align u &=\int_1^{\infty}\chi(t) e^{\xi(u) t}dt \leq
B_\epsilon\int_1^{\infty}\left(\frac{e^{\xi(u)+\epsilon/16}}{h(t)}\right)^tdt\\
&=
B_\epsilon\int_1^{\log(u)}\left(\frac{e^{\xi(u)+\epsilon/16}}{h(t)}\right)^tdt+B_\epsilon\int_{\log(u)}^{\infty}
\left(\frac{e^{\xi(u)+\epsilon/16}}{h(t)}\right)^tdt.\tag{28}\endalign$$
Considering the cases $t_0\leq \log(u)$ and $t_0>\log(u)$ (in
which case $q(t)$ is increasing on $[1,\log(u)]$), and using (27)
and our assumption on $\xi(u)$ we get that
$$\int_1^{\log(u)}\left(\frac{e^{\xi(u)+\epsilon/16}}{h(t)}\right)^tdt
\leq
\max(\log(u)e^{o(\log(u))},\log(u)\exp(-\epsilon/2\log(u)))=u^{o(1)},$$using
this, (28) and the assumption on $\xi(u)$ we deduce that
$$ u \leq B_\epsilon u^{o(1)} + B_\epsilon\int_{\log(u)}^{\infty} e^{-\epsilon/2 t}dt = u^{o(1)} +o(1),$$
which contradicts our hypothesis.

\noindent Now in both cases $h'(t)/h(t)\leq c/t$ for some positive
constant $c$ and for all $t$. Then integrating both sides gives
$h(t)\ll t^c$, and this with our result implies $\xi(u)\ll
\log(\log(u))$. Thus by Lemma 6.1 and Theorem 3 the Proposition
follows.

\enddemo

\enddocument